\documentclass[a4paper,12pt]{article}

\usepackage{amsthm,amsmath,amssymb}

\usepackage{graphicx}
\usepackage[colorlinks=true,citecolor=black,linkcolor=black,urlcolor=blue]{hyperref}

\theoremstyle{plain}
\newtheorem{theorem}{Theorem}[section]

\newtheorem{proposition}[theorem]{Proposition}

\theoremstyle{definition} 

\theoremstyle{definition}
\newtheorem{definition}[theorem]{Definition}
\newtheorem{example}[theorem]{Example}

\theoremstyle{remark}

\title{\bf Prime Number Sums}

\author{Paul Bradley\\
\small Senior Consultant\\[-0.8ex]
\small Jacobs\\[-0.8ex] 
\small Birmingham.\\
\small\tt paulmbradley82@gmail.com\\
}

\date{Jun 26, 2017}

\begin{document}

\maketitle

\section{Introduction}
There is a competition in which five contestants compete against each other over a series of rounds. Each of the five contestants is awarded a score between $1$ and $5$ points at the end of each round. Their score is based on the position they had finished in that round. After each round first place receives $5$ points, second receives $4$ points and so on down to last place receiving $1$ point.

After the first two rounds the scores looked like this:

\begin{table}[htbp]
\begin{center}
\begin{tabular}{@{\extracolsep{\fill}}| l | c | c | c | c | c | }
\hline
Round One       & $1$ & $2$  & $3$ & $4$ & $5$ \\ \hline
Round Two       & $3$ & $2$  & $1$ & $5$ & $4$ \\ \hline
Total           & $4$ & $4$  & $4$ & $9$ & $9$ \\ \hline
\end{tabular}
\caption{Scores after the first two rounds.}
\label{Tab1}
\end{center}
\end{table}

An interesting feature of this edition of the competition is that the totals are all composite numbers. Is it possible to have all of the totals equal to prime numbers after two rounds in much the same way?

After a little effort it will become apparent that it is possible to achieve such a result (see Table~\ref{Tab2}) and in fact there is only one way in which we can do this.

\begin{table}[htbp]
\begin{center}
\begin{tabular}{@{\extracolsep{\fill}}| l | c | c | c | c | c | }
\hline
Round One       & $1$ & $2$  & $3$ & $4$ & $5$ \\ \hline
Round Two       & $1$ & $5$  & $4$ & $3$ & $2$ \\ \hline
Total           & $2$ & $7$  & $7$ & $7$ & $7$ \\ \hline
\end{tabular}
\caption{All prime scores after the first two rounds.}
\label{Tab2}
\end{center}
\end{table}

This is not some special feature of $5$ contestants. Checking the cases for $1, 2, 3$ and $4$, it is possible to come up with a sequence of scores for the two rounds such that all contestants have a prime number total.

The next question is then obvious. If this is possible for up to $5$ contestants what is the smallest number of contestants for which this is not possible? Having found solutions which were unique for a number of the cases tested it appears that this is likely to fall down at some point.

\section{The Question}

The first thing to do is phrase the question mathematically.

\begin{definition}
Let $\mathbb{N}_n$ be the finite set of integers $\{1,2,...,n\}.$
\end{definition}

\begin{definition}
A permutation on the set $\mathbb{N}_n$ is a ordering of the numbers. Each number must appear exactly once in the ordering.
\end{definition}

\begin{example}
Let $n=3.$ Call the permutation $\pi$ and set $\pi = (1,3,2).$ Then this will send $1$ to the third position, $3$ to the second position and $2$ to the first position. So the sequence $ [ 1,2,3 ] $ would become $[2, 3, 1]$ after $\pi$ was applied. We write $\pi([1,2,3]) = [2,3,1]$ to apply the permutation to the entire sequence or $\pi(1)=3$ to apply the permutation to one number at a time.
\end{example}

Hence our question can be phrased as; 

\begin{center}
    
``Is there a permutation $\pi$ of the set $\mathbb{N}_n$ such that $k+\pi(k)$ is prime for each $1\le k\le n?$"
\end{center}

However, if there is no such permutation then we wish to consider the question of which is the smallest integer $n$ for which there is no such possible permutation.

We begin by considering the first few values of $n$. For $n=1$ this is trivial as $1+1=2.$

For the next few values of $n$ we have found solutions and presented them in Table $3$.

\begin{table}[htbp]
\begin{center}
\begin{tabular}{@{\extracolsep{\fill}}| l | c | c | c | c | } \hline

Round $1$       & $1$ & $2$ & &   \\ \hline
Round $2$       & $2$ & $1$  & & \\ \hline
Total     & $3$ & $3$  & &  \\ \hline
& & & & \\
\hline
Round $1$       & $1$ & $2$  & $3$ & \\ \hline
Round $2$       & $1$ & $3$  & $2$ & \\ \hline
Total     & $2$ & $5$  & $5$ & \\ \hline
& & & & \\
\hline
Round $1$       & $1$ & $2$  & $3$ & $4$  \\ \hline
Round $2$      & $2$ & $1$  & $4$ & $3$  \\ \hline
Total     & $3$ & $3$  & $7$ & $7$  \\ \hline
\end{tabular}
\caption{$n=2,3,4$.}
\label{Tab3}
\end{center}
\end{table}
Relating this problem to group theory allows us to state it in terms of permutation groups and to make use of the computational power of MAGMA.

Let $G \cong S(n)$ the full symmetric group on $n$ points. Let $G$ act naturally on $\mathbb{N}_n$. Then we aim to find $g \in G$ such that $i+i^g$ is prime for all $i \in \mathbb{N}_n.$

This rephrasing makes the MAGMA procedure below more clear.

\begin{example}

We can extend our search for the first value of $n$ where there is no such group element by writing a short MAGMA procedure.

\begin{verbatim}

AddPrimeTest:=procedure(n);
Z:=Integers();
G:=Sym(n);
for g in G do;
    S:=[];
	for t:=1 to n do;
		if IsPrime(Z!(t+t^g)) then;
			Append(~S, Z!(t+t^g));
		end if;
	end for;
    if #S eq n then;
        print g;
        break;
    end if;

end for;
end procedure;
\end{verbatim}
\end{example}

We can then quickly find solutions for the first few values of $n$. The above procedure was written to stop once it had found a solution and then print the appropriate permutation. The break command prevents the computer from finding any further solutions and saves memory and time.

Computation provides us with solutions for the first twenty values of $n$. This begs the question then, is there a value for which the task is impossible? There does not seem to be an obvious reason as to why there should be such a permutation in every case.

This leads us to our result.

\begin{proposition}
For all $n \in \mathbb{N}$ there exists a permutation $\pi_n$ of the set $\mathbb{N}_n$ such that $k+\pi_n(k)$ is prime for each $1\le k\le n.$
\end{proposition}

\begin{proof}
We proceed via strong induction. For the case $n=1$ letting $\pi_1$ be the identity permutation gives us the solution. For the case $n=2$ if we let $\pi_2 \in S(2)$ be $(1,2)$ then we have final scores of $1+2=2+1=3.$ For the case $n=3$ we can let $\pi_3=(2,3)$ which gives totals of $2,5$ and $5$. We now assume that the statement holds for all $n\le k-1$ for some $k>1$, that is for all $1\le i \le k-1$ there exists $\pi_i$. So we now consider the case for $n=k$. Here we choose the lowest prime $p$ greater than $k$. Assume now that $p=k+a$ where $a<k$. Then we can form the permutation $\pi_a^{k}$ which fixes all $1\le i  <a $ and pairs the remaining points as shown.

\begin{table}[htbp]
\begin{center}
\begin{tabular}{@{\extracolsep{\fill}}| l | c | c | c | c | c | }
\hline
Round $1$           & $a$ & $a+1$  & $...$ & $k-1$ & $k$ \\ \hline
Round $2$           & $k$ & $k-1$  & $...$ & $a+1$ & $a$ \\ \hline
Total               & $k+a$ & $k+a$  & $k+a$ & $k+a$ & $k+a$ \\ \hline
\end{tabular}
\label{Tab6}
\end{center}
\end{table}

The fixed points $\{1,2,...,a-1\}$, if they exist, can now be solved by using $\pi_{a-1}$ which exists by assumption. Combining these two permutations by performing $\pi_a^{k}$ followed by $\pi_{a-1}$ gives us  $\pi_k$ and hence our result. We now only need to consider the case where $p\ge 2k$. However this case is ruled out by the Bertrand–-Chebyshev Theorem~\cite{Cheby} which states that for every $n>1,$ there is some prime $p$ with $n<p<2n$.
\end{proof}

Given the construction in our proof it looks like the number of solutions to this problem for each $n$ is related to the number of primes less than $2n$ and the number of solutions for the remaining points after pairing to make the primes below $2n$. Removing the break command in the procedure allows us to test for the number of solutions for small $n$. The results are collated in Table~\ref{Tab7}.

\begin{table}[htpb]
\begin{center}
\begin{tabular}{@{\extracolsep{\fill}}| l | c | c | }
\hline
$n$         & Number of Permutations &  Number of Primes $<2n$\\ \hline
$1$         & $1$ &  $0$ \\ \hline
$2$         & $1$ &  $2$ \\ \hline
$3$         & $1$ &  $3$ \\ \hline
$4$         & $4$ &  $4$ \\ \hline
$5$         & $1$ &  $4$ \\ \hline
$6$         & $9$ &  $5$ \\ \hline
$7$         & $4$ &  $6$ \\ \hline
$8$         & $36$ &  $6$ \\ \hline
$9$         & $36$ &  $7$ \\ \hline
$10$         & $676$ &  $8$ \\ \hline
\end{tabular}
\caption{Number of distinct permutation solutions for $n\le 10$.}
\label{Tab7}
\end{center}
\end{table}

Further details of the sequence for number of solutions in each case can be found in \cite{Seq}.

\end{document}